\documentclass[12pt]{article}
\usepackage{amsmath,amsthm,latexsym,amssymb,pgf} 
\usepackage{amsmath,amsthm,latexsym,amssymb,pgf,makeidx,
dsfont,authblk,boxedminipage,tikz}
\usetikzlibrary{arrows,decorations.pathmorphing,backgrounds,positioning,fit,petri,calc}
\usepackage{latexsym,amsmath} 
\begin{document}
\newtheorem{defi}{Definition}
\newtheorem{theo}{Theorem}
\newtheorem{coro}{Corollary}
\newtheorem{prop}{Proposition}
\newtheorem{obs}{Observation}
\newtheorem{rem}{Rem:}
\newtheorem{la}{Lemma}
\newtheorem{con}{Conjecture}
\newtheorem{problem}{Problem}
\newtheorem{exa}{Example}
\begin{center}
{\Large {\bf On the Wiener Index of Orientations of Graphs}} \\[2mm]
{\large Peter Dankelmann\footnote{Financial support by the South African
National Research Foundation, grant 118521, is gratefully acknowledged}
\\ University of Johannesburg} \\
\end{center}

\begin{abstract}
The Wiener index of a strong digraph $D$ is defined as the sum of the distances 
between all ordered pairs of vertices. This definition has been extended to
digraphs that are not necessarily strong by defining the distance from a vertex
$a$ to a vertex $b$ as $0$ if there is no path from $a$ to $b$ in $D$. 

Knor, \u{S}krekovski and Tepeh [Some remarks on Wiener index of oriented graphs.
Appl.\ Math.\ Comput.\ {\bf 273}] considered orientations of graphs with maximum 
Wiener index. The authors conjectured that for a given tree $T$, an orientation 
$D$ of $T$ of maximum Wiener index always contains a vertex $v$ such that 
for every vertex $u$, there  is either a $(u,v)$-path or a $(v,u)$-path in $D$.
In this paper we disprove the conjecture. 

We also show that the problem of
finding an orientation of maximum Wiener index of a given graph is NP-complete,
thus answering a question by Knor, \u{S}krekovski and Tepeh 
[Orientations of graphs with maximum Wiener index. Discrete Appl.\ Math.\  211]. 

We briefly discuss the corresponding problem of 
finding an orientation of minimum Wiener index of a given graph, and show that
the special case of deciding if a given graph on $m$ edges has an orientation of
Wiener index $m$ can be solved in time quadratic in $n$. 
\end{abstract}


\section{Introduction}

The Wiener index $W(G)$ of a connected graph $G$ is defined as the sum of the distances
between all unordered pairs of vertices, i.e.,
\[ W(G) = \sum_{\{u,v\} \subseteq V} d_G(u,v), \]
where $V$ is the vertex set of $G$ and $d_G(u,v)$ denotes the distance between 
$u$ and $v$, i.e., the length a shortest $(u,v)$-path. The Wiener index, originally
introduced as a tool in chemistry \cite{Wie1947}, has been studied extensively in the 
mathematical and chemical literature and is arguably one of the most applicable
graph invariants. 

The definition of Wiener index extends naturally to strong digraphs. If $D$ is a 
strong digraph with vertex set $V$, then the Wiener index of $D$ is defined as 
\[ W(D) = \sum_{(u,v) \in V\times V} d_D(u,v), \]
where $V \times V$ is the set of all ordered pairs of vertices of $D$. 
The Wiener index of a graph (digraph) is closely related to the average distance, defined as the 
arithmetic mean of the distances between all unordered (ordered) pairs of distinct vertices.
For recent results on the Wiener index of strong digraphs see, for example,  
\cite{HuaLin2016} and \cite{Dan2019}.

Knor, \u{S}krekovski and Tepeh \cite{KnoSkrTep2016} extended the definition of 
the Wiener index to digraphs that are not necessarily strong by defining
the distance $d_D(u,v)$ as the length of a shortest $(u,v)$-path if 
$D$ contains a $(u,v)$-path, and as $0$ if no $(u,v)$-path exists in $D$. 
In \cite{KnoSkrTep2016-2} the same authors gave further results on the Wiener 
index, for example they showed that the Wiener index of a  tournaments of
order $n$ cannot exceed $\binom{n+1}{3}-1$, which had previously been proved
by Plesn\' \i k \cite{Ple1984} for tournaments that are strong. 

The Wiener index of orientations of graphs was investigated first in \cite{DanOelWu2004}, 
where strong orientations  of a given graph that minimise the Wiener index were considered. 
Applying their more general definition of distance in digraphs,  
Knor, \u{S}krekovski and Tepeh \cite{KnoSkrTep2016, KnoSkrTep2016-2} considered 
(not necessarily strong) orientations of graphs that maximise or minimise the Wiener index. 
Among several other results, they demonstrated that it is not true in general that 
an orientation of a given graph maximising the Wiener index is necessarily strong. 

In this paper we resolve two open questions related to orientations of maximum Wiener index  
by Knor, \u{S}krekovski and Tepeh. Considering trees, they conjectured the following result.

\begin{con}[\cite{KnoSkrTep2016}]  \label{conj:core-vertex}
Let $T$ be a tree. If $D$ is an orientation of $T$ that maximises $W(D)$, then
there exists a vertex $v$ in $D$ such that for every vertex $u$ there exists
a $(u,v)$-path or a $(v,u)$-path in $D$. 
\end{con}

We prove that this conjecture is not true in general. 
We also consider the question whether finding an orientation of maximum Wiener index of 
a given graph is NP-hard, posed in \cite{KnoSkrTep2016-2}, and answer it  
in the affirmative. 
Finally, we briefly discuss the problem of finding an orientation of minimum Wiener
index of a given graph.

We use the following notation. 
We denote the vertex set of a graph (digraph, mixed graph) $G$ by $V(G)$, and the 
edge set or arc set by $E(G)$, leaving out the argument $G$ if there is no danger
of confusion. Generally, $uv$ denotes an undirected edge, while 
$\overrightarrow{uv}$ denotes a directed edge which is directed from $u$ to $v$.  

The converse of a digraph $D$ is the digraph obtained from $D$ by reversing the 
direction of every arc of $D$. 

A Hamiltonian path in a graph $G$ is a path that contains all vertices of $G$.

\section{A Counter-example to Conjecture \ref{conj:core-vertex} }

Conjecture \ref{conj:core-vertex} is supported by results in \cite{KnoSkrTep2016}, 
which show that it is true for some subclasses of trees. 
Since it is reasonable to expect that an orientation of a tree maximising the Wiener 
index also maximises the number of pairs of vertices $(u,v)$ between which there exists
a path, the following result due to Henning and Oellermann \cite{HenOel2004} which gives 
further support to Conjecture \ref{conj:core-vertex}.

\begin{theo}[\cite{HenOel2004}] 
Let $T$ be a tree, and $D$ an orientation of $T$ that maximises the number 
of ordered pairs $(u,v)$ of vertices of $D$ for which there exists a $(u,v)$-path
in $D$. Then $D$ contains a vertex $w$ so that for every vertex $u$ there exists
a $(u,w)$-path or a $(w,u)$-path in $D$. 
\end{theo}

Nevertheless, we found that Conjecture \ref{conj:core-vertex} is not true
in general. In this section we present an infinite family of counter-examples
to  Conjecture \ref{conj:core-vertex}.

Let $D$ be an orientation of a tree $T$. Following \cite{KnoSkrTep2016}, we say that
$D$ is zig-zag if $T$ contains a path $P$ whose edges change their direction in $D$ at least
twice as $P$ is traversed. It was observed in \cite{KnoSkrTep2016} that an orientation
$D$ of a tree is not zig-zag if and only if $D$ has a vertex $w$ so that for every 
vertex $u$ of $D$ there exists either a $(u,w)$-path or a $(w,u)$-path. 
Conjecture \ref{conj:core-vertex} was given in \cite{KnoSkrTep2016} in an equivalent 
form, stating that for a given tree $T$, every orientation of $T$ that
has maximum Wiener index is not zig-zag.

We construct a family of counter-examples to Conjecture \ref{conj:core-vertex} as follows. 
Let $k \in \mathbb{N}$ be a multiple of $3$. Let $T_k$ be the tree obtained from a path
of order $k$ with vertices $w_1,w_2,\ldots,w_k$ by 
appending vertices $u_1, u_2,\ldots, u_{k^2/9}$ to $w_1$, 
appending a path $x_1,x_2,x_3,x_4,x_5$ to $w_2$, and a single vertex $y_1$ to $w_3$.
A sketch of the tree $T_k$ is shown in Figure \ref{fig:tree-Tk}.

  \begin{figure}[h]
  \begin{center}
\begin{tikzpicture}
  [scale=0.6,inner sep=1mm, 
   vertex/.style={circle,thick,draw}, 
   thickedge/.style={line width=2pt}] 
    \node[vertex] (a1) at (-1,0) [fill=white] {};
    \node[vertex] (b1) at (2,2) [fill=white] {};    
    \node[vertex] (b2) at (4,2) [fill=white] {};
    \node[vertex] (b3) at (6,2) [fill=white] {};
    \node[vertex] (b4) at (8,2) [fill=white] {};
    \node[vertex] (b5) at (10,2) [fill=white] {};
    \node[vertex] (b6) at (14,2) [fill=white] {};    
    \node[vertex] (c1) at (4,4) [fill=white] {};
    \node[vertex] (c2) at (6,4) [fill=white] {};    
    \node[vertex] (d1) at (4,6) [fill=white] {};
    \node[vertex] (e1) at (4,8) [fill=white] {};
    \node[vertex] (f1) at (-1,10) [fill=white] {};   
    \node[vertex] (f2) at (4,10) [fill=white] {};       
    \node[vertex] (g1) at (-1,12) [fill=white] {};    
    \node[vertex] (g2) at (4,12) [fill=white] {};        
    \draw[black, thick] (b1)--(b2)--(b3)--(b4)--(b5) 
        (b2)--(c1)--(d1)--(e1)--(f2)--(g2)
        (a1)--(b1)  (f1)--(b1)   (g1)--(b1)  
        (b3)--(c2);
   \draw[black, thick, dashed] (b5)--(b6);
   \node[left] at (-1.2,0) {$u_{k^2/9}$};
   \node[left] at (-1.2,10) {$u_{2}$};
   \node[left] at (-1.2,12) {$u_{1}$};   
   \node[left] at (3.8,4) {$x_1$};      
   \node[left] at (3.8,6) {$x_2$};   
   \node[left] at (3.8,8) {$x_3$};   
   \node[left] at (3.8,10) {$x_4$};   
   \node[left] at (3.8,12) {$x_5$};   
   \node[right] at (6.2,4) {$y_1$};     
   \node[below] at (2,1.8) {$w_1$};     
   \node[below] at (4,1.8) {$w_2$};  
   \node[below] at (6,1.8) {$w_3$};  
   \node[below] at (8,1.8) {$w_4$};         
   \node[below] at (10,1.8) {$w_5$};    
   \node[below] at (14,1.8) {$w_k$};    
   \draw[black, dashed] (-1,2)--(-1,8);            
\end{tikzpicture}
\end{center}
\caption{The tree $T_k$.}
\label{fig:tree-Tk}
\end{figure}
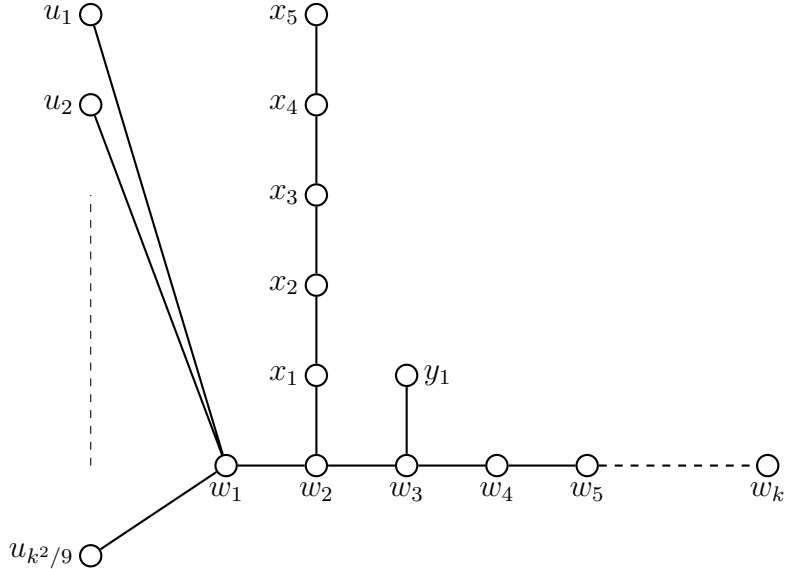

  \begin{figure}[h]
  \begin{center}
\begin{tikzpicture}
  [scale=0.6,inner sep=1mm, 
   vertex/.style={circle,thick,draw}, 
   thickedge/.style={line width=2pt}] 
    \begin{scope}[>=triangle 45]
       
    \node[vertex] (a1) at (-1,0) [fill=white] {};
    \node[vertex] (b1) at (2,2) [fill=white] {};    
    \node[vertex] (b2) at (4,2) [fill=white] {};
    \node[vertex] (b3) at (6,2) [fill=white] {};
    \node[vertex] (b4) at (8,2) [fill=white] {};
    \node[vertex] (b5) at (10,2) [fill=white] {};
    \node[vertex] (b6) at (14,2) [fill=white] {};    
    \node[vertex] (c1) at (4,4) [fill=white] {};
    \node[vertex] (c2) at (6,4) [fill=white] {};    
    \node[vertex] (d1) at (4,6) [fill=white] {};
    \node[vertex] (e1) at (4,8) [fill=white] {};
    \node[vertex] (f1) at (-1,10) [fill=white] {};   
    \node[vertex] (f2) at (4,10) [fill=white] {};       
    \node[vertex] (g1) at (-1,12) [fill=white] {};    
    \node[vertex] (g2) at (4,12) [fill=white] {};        
    \draw[black, thick, ->] (a1)--(b1);
    \draw[black, thick, ->] (f1)--(b1);   
    \draw[black, thick, ->](g1)--(b1);
   \draw[black, thick, ->] (b1)--(b2);   
   \draw[black, thick, ->] (b2)--(b3);   
   \draw[black, thick, ->] (b3)--(b4);   
   \draw[black, thick, ->] (b4)--(b5);  
   \draw[black, thick, ->] (b2)--(c1);  
   \draw[black, thick, ->] (c1)--(d1);   
   \draw[black, thick, ->] (d1)--(e1);   
   \draw[black, thick, ->] (e1)--(f2);      
   \draw[black, thick, ->] (f2)--(g2);                   
   \draw[black, thick, dashed, ->] (b5)--(b6);                   
   \draw[black, thick, ->] (c2)--(b3);
   \node[left] at (-1.2,0) {$u_{k^2/9}$};
   \node[left] at (-1.2,10) {$u_{2}$};
   \node[left] at (-1.2,12) {$u_{1}$};   
   \node[left] at (3.8,4) {$x_1$};      
   \node[left] at (3.8,6) {$x_2$};   
   \node[left] at (3.8,8) {$x_3$};   
   \node[left] at (3.8,10) {$x_4$};   
   \node[left] at (3.8,12) {$x_5$};   
   \node[right] at (6.2,4) {$y_1$};     
   \node[below] at (2,1.8) {$w_1$};     
   \node[below] at (4,1.8) {$w_2$};  
   \node[below] at (6,1.8) {$w_3$};  
   \node[below] at (8,1.8) {$w_4$};         
   \node[below] at (10,1.8) {$w_5$};    
   \node[below] at (14,1.8) {$w_k$};    
   \draw[black, dotted] (-1,2)--(-1,8);            
   \end{scope}
\end{tikzpicture}
\end{center}
\caption{The orientation of $D_k$ of $T_k$.}
\label{fig:optimal-orientation}
\end{figure}
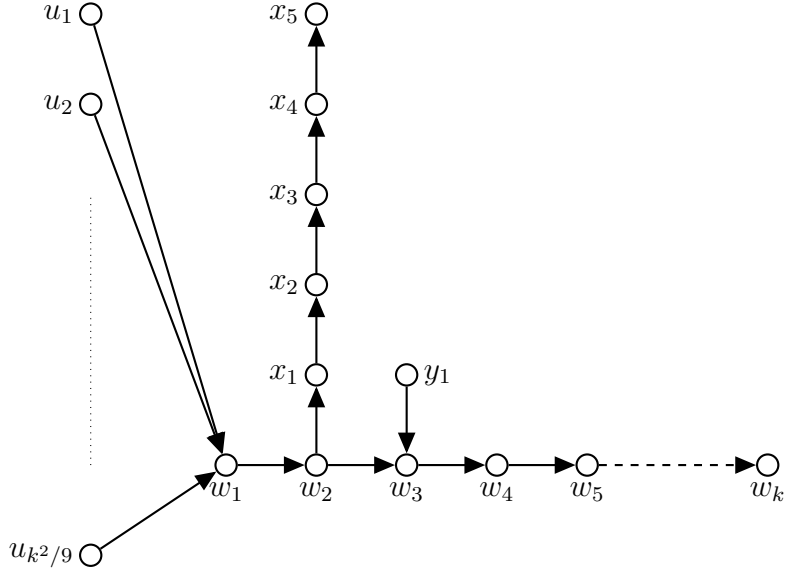

Let $D_k$ be the orientation of $T_k$  shown in Figure \ref{fig:optimal-orientation}, i.e., 
the edges of the path $w_1, w_2,\ldots, w_k$ are
oriented towards $w_k$, each edge $u_iw_1$ is oriented towards $w_1$, the edges of the path
$w_2, x_1, x_2,\ldots,x_5$ are oriented towards $x_5$, and the edge $y_1w_3$ is oriented towards $w_3$. 
Notice that the edges of the $(x_5,y_1)$-path change their direction twice as the path
is traversed, hence $D_k$ is zig-zag.

It will be convenient to consider mixed graphs, which are 
a common generalisation of graphs and digraphs. 
A mixed graph consists of a set of vertices and a set of edges, where each edge may or 
may not have a direction.  For vertices $u,v$ of a mixed graph $G$, a $(u,v)$-path is
a sequence $v_0, v_1,\ldots,v_k$ of vertices with $u=v_0$ and $v=v_k$ so that for 
each $i \in \{0,1,\ldots,k-1\}$, the directed edge $\overrightarrow{v_iv_{i+1}}$ 
or the undirected edge $v_iv_{i+1}$ is in $G$. As usual, the distance from $u$ to 
$v$ is the minimum number of edges on a $(u,v)$-path, which we denote by $d_G(u,v)$. 
By a partial orientation of a mixed graph $G$ we mean a mixed graph obtained
from $G$ by orienting some of the undirected edges of $G$.

In order to compare the Wiener index of a tree and its (partial) orientations, 
we introduce the following modification of the Wiener index
For a mixed graph $G$ define 
\[ W^{max}(G) := \sum_{ \{u,v\} \subseteq V } \max\{ d_G(u,v), d_G(v,u) \}. \] 
For disjoint subsets $A,B \subseteq V$  
we use the notation $W^{max}_G(A,B) := \sum_{\ a\in A, b \in B} \max\{d_G(a,b), d_G(b,a)\}$.

\begin{la} \label{la:Wmax}
(a) If $T$ is a tree, then 
\[ W^{max}(T) = W(T). \]
(b) If $D$ is an orientation of a tree $T$, then 
\[ W^{max}(D) = W(D). \] 
(c) If $D_1$ is a partial orientation of z tree $T$, and $D_2$ a partial orientation of $D_1$,
then 
\[ W^{max}(D_2)  \leq W^{max}(D_1) \leq   W(T). \]
(d) Let $T$ be a tree, $D_1$ a partial orientation of $T$, and $D_2$ a partial orientation of $D_1$.
If $A$ and $B$ are disjoint sets of vertices of $T$ such that there is no path in $D_2$ 
between a vertex in $A$ and a vertex in $B$ in either direction,  then 
\[ W^{max}(D_2) \leq W^{max}(D_1) - W^{max}_{D_1}(A,B). \]
\end{la}

{\bf Proof:}
We denote the vertex set of $T$ and its  (partial) orientations by $V$. \\[1mm]
(a) Since $T$ is an undirected graph, we have $d_T(u,v) = d_T(v,u) = \max\{ d_T(u,v), d_T(v,u)\}$
for any two vertices $u, v$ of $T$. Summation over all subsets $\{u,v\} \subseteq V$
yields that $W^{max}(T) = W(T)$. \\[1mm]
(b) If $u,v \in V$ are two vertices of $T$, then the $(u,v)$-path and the 
$(v,u)$-path are unique in $T$. At most one these two paths is also a path in $D$, 
thus $d_D(u,v)=0$ or $d_D(v,u)=0$. Hence 
$\max\{d_D(u,v), d_D(v,u)\} = d_D(u,v)+d_D(v,u)$. Summation over all subsets
$\{u,v\} \subseteq V$ yields the statement of (b). \\[1mm]
(c) Let $u,v$ be two vertices of $T$. If there is a $(u,v)$-path in $D_2$, then this path 
is also a $(u,v)$-path in $D_1$, and if there is a $(u,v)$-path in $D_1$, then this path is 
also a $(u,v)$-path in $T$. In all cases, the $(u,v)$-path is unique if it exists. It follows that
$d_{D_2}(u,v) \leq d_{D_1}(u,v) \leq d_T(u,v)$. Hence, for all $u,v \in V$, 
\begin{equation}  \label{eq:two-sub-orientations}
\max\{d_{D_2}(u,v), d_{D_2}(v,u)\} \leq \max\{d_{D_1}(u,v), d_{D_1}(v,u)\}
                  \leq d_{T}(u,v). 
\end{equation}                  
Summation over all subsets $\{u,v\} \subseteq V(T)$ yields (c).  \\[1mm]              
(d) Let ${\cal C}$ be the set of all $2$-vertex subsets of $V$ which do not consist of a 
vertex of $A$ and a vertex of $B$. Then $\max\{ d_{D_2}(u,v), d_{D_2}(v,u) \}=0$ 
whenever $\{u,v\} \notin {\cal C}$. Hence, by \eqref{eq:two-sub-orientations},  
\begin{eqnarray*} W^{max}(D_2) 
  & = &  \sum_{ \{u,v\} \in {\cal C} } \max\{ d_{D_2}(u,v), d_{D_2}(v,u) \} \\
  & \leq &  \sum_{ \{u,v\} \in {\cal C} } \max\{ d_{D_1}(u,v), d_{D_1}(v,u) \} \\
  &  =  & W^{max}(D_1) - W^{max}_{D_1}(A,B), 
\end{eqnarray*}  
as desired. \hfill $\Box$ \\

\begin{theo} \label{theo:counter-example-to-Conjecture1} 
Let $k \in \mathbb{N}$ be a multiple of $3$. 
Let $D_k$ be the orientation of $T_k$ as defined above, and let $D$ be any orientation of $T_k$. 
If $k$ is sufficiently large, then  
\begin{equation} \label{eq:Dk-is-maximal} 
W(D) \leq W(D_k), 
\end{equation}
with equality if and only if $D$ equals $D_k$ or the converse of $D_k$. 
\end{theo}

{\bf Proof:} 
We first determine the Wiener indices of $T_k$ and $D_k$. 
Tedious but straightforward calculations show that 
\begin{equation} \label{eq:value-of-W(Tk)}  
W(T_k)  = \frac{11}{162}k^4 + \frac{2}{9}k^3 + \frac{55}{9}k^2 + \frac{35}{6}k + 61 
\end{equation}
and 
\begin{equation} \label{eq:value-of-W(Dk)}
W(D_k) = \frac{1}{18}k^4 + \frac{2}{9}k^3 + \frac{59}{18}k^2 - \frac{5}{3}k + 56. 
\end{equation}
It suffices to prove the theorem for orientations of $T_k$ of maximum Wiener index. 
Let $D$ be such an orientation of $T_k$. 
Then $W^{max}(D)$ is maximum among all orientations of $T_k$ by Lemma \ref{la:Wmax}(b). 
We may further assume that $\overrightarrow{w_1w_2} \in E(D)$; otherwise we consider the 
converse of $D$. We prove that $D=D_k$ for sufficiently large values of $k$.

Let $V$ be the common vertex set of $T_k$, $D_k$ and $D$. 
We partition $V$ into four sets, $U$, $W$, $X$ and $Y$, where 
$U=\{u_1,u_2,\ldots,u_{k^2/9}\}$, $W=\{w_1, w_2,\ldots,w_k\}$,
$X=\{x_1, x_2, x_3, x_4, x_5\}$, and $Y=\{y_1\}$. \\[1mm]
{\sc Claim 1:} $D$ contains a $(w_1, w_4)$-path. \\
Suppose to the contrary that $D$ contains no $(w_1, w_4)$-path. Since 
$\overrightarrow{w_1w_2} \in E(D)$, there is also $(w_4,w_1)$-path in $D$. 
Hence $D$ contains no path between a vertex in $U$ and a vertex in $W-\{w_1, w_2, w_3\}$
in either direction. Applying Lemma \ref{la:Wmax}(c) to the two suborientations
$T_k$ and $D$ of $T_k$ we obtain that 
\[ W^{max}(D) \leq W^{max}(T_k) - W^{max}_{T_k}(U,W-\{w_1, w_2,w_3\}). \]    
An easy calculation shows that 
$W^{max}_{T_k}(U,W-\{w_1, w_2,w_3\}) = \frac{1}{18}k^4 + \frac{1}{18}k^3 - \frac{2}{3}k^2$. 
Since $W^{max}(T_k)=W(T_k)$ by Lemma \ref{la:Wmax}(a), we thus obtain from  
\eqref{eq:value-of-W(Tk)} that 
\begin{eqnarray*} 
W^{max}(D) & \leq & 
   \big(  \frac{11}{162}k^4 + \frac{2}{9}k^3 + \frac{55}{9}k^2 + \frac{35}{6}k + 61  \big) 
               - \big(  \frac{1}{18}k^4 + \frac{1}{18}k^3 - \frac{2}{3}k^2 \big) \\
   & = & \frac{2}{162}k^4 + \frac{1}{6}k^3 + \frac{61}{9}k^2 + \frac{35}{6}k + 61. 
\end{eqnarray*}             
Comparing this with the right hand side of \eqref{eq:value-of-W(Dk)}, it is easy to see  that 
$W^{max}(D) < W^{max}(D_k)$ for sufficiently large $k$. This contradiction to 
the maximality of $W^{max}(D)$ proves Claim 1. \\
{\sc Claim 2:} $D$ contains a $(w_1,w_k)$-path. \\
Suppose to the contrary that
$D$ does not contain a $(w_1,w_k)$-path. Let $i$ be the smallest value for which 
there exists no $(w_1,w_i)$-path. By Claim 1 we have $i \geq 5$. Then $w_{i-1}$ is adjacent from 
$w_{i-2}$ and $w_{i}$, and not adjacent to any vertex in $D$. Reversing all arcs
along the path $w_{i-1}, w_{i},\ldots,w_k$ does not reduce $\max\{d_D(u,v), d_D(v,u)\}$ for any
pair of vertices, but increases $\max\{d_D(w_{i-2}, w_i), d_D(w_i, w_{i-2}) \}$, 
contradicting the maximality of $W^{max}(D)$. Claim 2 follows. \\[1mm]
{\sc Claim 3:} $\overrightarrow{u_iw_1} \in E(D)$ for all $i \in \{1,2,\ldots,\frac{k^2}{9}\}$. \\ 
Suppose to the contrary that $\overrightarrow{w_1u_i} \in E(D)$ for some $i$.
Reversing the arc $\overrightarrow{w_1u_i}$ creates 
paths from $u_i$ to all vertices of $W$, whose total length is  
$\sum_{j=1}^k d_{T_k}(u_i,w_j) = \frac{k(k+1)}{2}$,
but destroys only paths from vertices in $\{w_1\} \cup (U-\{u_i\})$ to $u_i$,
whose total length is not more than $1+2(|U|-1)= \frac{2}{9}k^2-1$.  
Since $ \frac{k(k+1)}{2} > 2\frac{k^2}{9}-1$, reversing  the 
arc $\overrightarrow{w_1u_i}$ increases the Wiener index, a contradiction
to the maximality of $W^{max}(D)$. Claim 3 follows. \\[1mm]
Let $D_k'$ be the partial orientation of $T_k$ in which for $i=1,2,\ldots,k-1$ the edge
$w_iw_{i+1}$ receives the orientation $\overrightarrow{w_iw_{i+1}}$, and 
for each $u_j \in U$ the edge $u_jw_1$ received the orientation $\overrightarrow{u_jw_1}$, 
while the remaining edges have not been oriented. It follows from Claims 2 and 3 that
$D$ is a suborientation of $D_k'$. A simple calculation shows that
\[ W^{max}(D_k') = \frac{1}{18}k^4 + \frac{2}{9}k^3 + \frac{56}{9}k^2 + \frac{35}{6}k + 61. \]
{\sc Claim 4:} $D$ contains a $(w_2,x_5)$-path. \\
We first show that $\overrightarrow{w_2x_1} \in E(D)$. Suppose to the contrary that
$\overrightarrow{x_1w_2} \in E(D)$. Then $D$ contains no path between a vertex of $X$ and 
a vertex of $U$ in either direction. Applying Lemma \ref{la:Wmax}(d) to $D_k'$ and $D$ yields
\begin{eqnarray*} 
W^{max}(D) & \leq & W^{max}(D_k') - W^{max}_{D_k'}(U,X)  \\
   & = & \big( \frac{1}{18}k^4 + \frac{2}{9}k^3 + \frac{56}{9}k^2 + \frac{35}{6}k + 61 \big)
     - \frac{25}{9} k^2 \\
   & = &  \frac{1}{18}k^4 + \frac{2}{9}k^3 + \frac{31}{9}k^2 + \frac{35}{6}k + 61,
\end{eqnarray*}    
and so, since $W^{max}(D_k) = \frac{1}{18}k^4 + \frac{2}{9}k^3 + \frac{59}{18} k^2 + O(k)$, it follows that 
$W^{max}(D) < W^{max}(D_k)$ for sufficiently large $k$, a contradiction to the maximality
of $W^{max}(D)$. Hence 
$\overrightarrow{w_2x_1} \in E(D)$. Similar arguments as in the proof of Claim 2 now
prove that $D$ contains a path from $w_2$ to $x_5$. Claim 4 follows. \\[1mm]
Let $D_k''$ be the partial orientation of $D_k'$ in which the edges of the $(w_2,x_1)$-path 
are oriented towards $x_5$. \\[1mm]
It follows from Claims 1 to 4 that $D$ is an orientation of $D_k''$. \\[1mm]
{\sc Claim 5:} $D=D_k$. \\
Only the edge $y_1w_3$ of $D_k''$ has not received an orientation. Hence $D_k''$ has two
orientations. $D_k$ (in which $y_1w_3$ receives the orientation $\overrightarrow{y_1w_3}$)
and the orientation in which $y_1w_3$ receives the orientation $\overrightarrow{w_3y_1}$, which
we denote by $D_k'''$.  Clearly, $d_{D_k}(u,v)=d_{D_k'''}(u,v)$ for all  
$u,v \in V-\{y_1\}$. Hence
\begin{eqnarray*} 
W^{max}(D_k) - W^{max}(D_k''') 
 & = & W^{max}_{D_k}(\{y_1\}, V-\{y_1\}) - W^{max}_{D_k'''}(\{y_1\}, V-\{y_1\}) \\
 & = & \sum_{v \in \{w_3, w_4,\ldots,w_k\} } d_{D_k}(y_1,v) 
            - \sum_{v \in U \cup X \cup \{w_1,w_2,w_3\} } d_{D_k'''}(v, y_1)  \\
 & = & (\frac{1}{2}k^2 - 2k + \frac{3}{2}) - (\frac{4}{9}k^2 - 6), 
\end{eqnarray*}
which is positive for $k$ sufficiently large. Hence $W^{max}(D_k) > W^{max}(D_k''')$. Claim 5 and 
thus the theorem follows. 
\hfill $\Box$

Since by Theorem \ref{theo:counter-example-to-Conjecture1} the only orientations of 
$T_k$ that maximise the Wiener index are $D_k$ and its converse, and both are zigzag,
it follows that Conjecture \ref{conj:core-vertex} is not true in general.

\section{Complexity of finding an orientation of maximum Wiener index}

In this section we consider the problem of finding a (not necessarily strong) 
orientation of a given graph that maximises the Wiener index. 
Knor, \u{S}krekovski and Tepe \cite{KnoSkrTep2016-2} asked whether this problem
is NP-hard, and we answer this question in the affirmative. 
Specifically, we consider the decision problem \\[2mm]
{\sc Wiener-Orientation:}  Given a graph $G$ and an integer $M$. Does $G$ have a
(not necessarily strong) orientation $D$ with $W(D) \geq M$?  \\[2mm]
We prove the NP-completeness of {\sc Wiener-Orientation} by a transformation
from the NP-complete problem {\sc Hamiltonian $(a,b)$-Path}, defined below. \\[2mm]
{\sc Hamiltonian $(a,b)$-path:} Given a graph $G$ and two vertices
$a$ and $b$ of $G$. Does $G$ have a Hamiltonian path that begins in $a$ and ends in $b$? \\[2mm]
In our proof we use the following notation. Let $D$ be a digraph with vertex set $V$
and let $A$ and $B$ be disjoint subsets of $V$. Then we write
$W_D(A,B)$ for $\sum_{a\in A, b \in B} d_D(a,b)$, and 
$W_D(A)$ for $\sum_{(a_1, a_2) \in A \times A} d_D(a_1, a_2)$.

\begin{figure}[h]
  \begin{center}
\begin{tikzpicture}
  [scale=0.6,inner sep=1mm, 
   vertex/.style={circle,thick,draw}, 
   thickedge/.style={line width=2pt}] 
    \node[vertex] (a1) at (0,0) [fill=white] {};
    \node[vertex] (a2) at (14,0) [fill=white] {};    
    \node[vertex] (b1) at (2,3) [fill=white] {};    
    \node[vertex] (b2) at (4,3) [fill=white] {};        
    \node[vertex] (b3) at (10,3) [fill=white] {};        
    \node[vertex] (b4) at (12,3) [fill=white] {}; 
    \node[vertex] (c1) at (0,5) [fill=white] {};
    \node[vertex] (c2) at (14,5) [fill=white] {};        
    \node[vertex] (d1) at (0,6) [fill=white] {};
    \node[vertex] (d2) at (14,6) [fill=white] {};            

    \draw[black, thick] (a1)--(b1)--(b2)  (b3)--(b4)--(a2); 
    \draw[black, thick] (c1)--(b1)  (b4)--(c2);   
    \draw[black, thick] (d1)--(b1)  (b4)--(d2);

   \node[left] at (-0.2,0) {$a_{n^3}$};
   \node[left] at (-0.2,5) {$a_{2}$};
   \node[left] at (-0.2,6) {$a_{1}$};   
   \node[below] at (2.08,2.8) {$a_{0}$};  
   \node[below] at (4,2.8) {$a$};        
   
   \node[right] at (14.2,0) {$b_{n^3}$};
   \node[right] at (14.2,5) {$b_{2}$};
   \node[right] at (14.2,6) {$b_{1}$};  
   \node[below] at (10,2.8) {$b$};  
   \node[below] at (12,2.8) {$b_0$};        
   
   \draw[black, rounded corners] (3.5,1) rectangle (10.5,5);
   \node[below] at (7,2.5) {$G$};      
    
   \draw[black, dotted] (0,0.7)--(0,4.3);    
   \draw[black, dotted] (14,0.7)--(14,4.3);            
\end{tikzpicture}
\end{center}
\caption{The graph $G^{a,b}$.}
\label{fig:graph-Gab}
\end{figure}
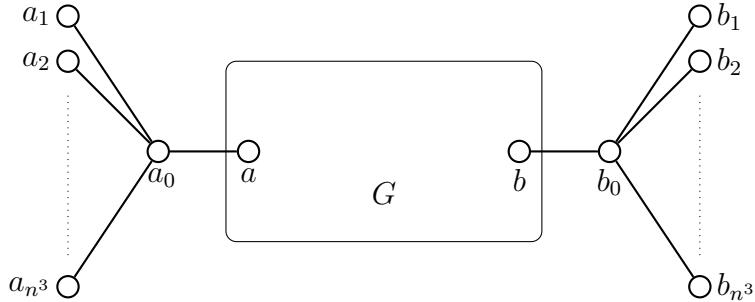

Given a graph $G$ of order $n$ and two vertices $a$ and $b$ of $G$, we define 
$G^{a,b}$ to be the graph of order $n^3+n+2$ obtained from $G$ by adding $2n^3+2$ 
new vertices $a_0, a_1,\ldots, a_{n^3}$ and $b_0, b_1, \ldots,b_{n^3}$ and
edges $aa_0$ and $bb_0$, as well as edges $a_0a_i$ and $b_0b_i$ for $i=1,2,\ldots,n^3$. 
A sketch of the graph $G^{a,b}$ is given in Figure  \ref{fig:graph-Gab}.

\begin{la} \label{la:Hamilton-if-large-W}
Let $M(n):= n^7 + 3n^6 + 2n^4 + 4n^3 +n+1$. Then there exists $N_0 \in \mathbb{N}$
such that for every graph $G$ of order $n$, where $n\geq N_0$, with vertices $a,b$, 
the following are equivalent: \\  
(i) $G$ has a Hamiltonian path from $a$ to $b$, \\
(ii) $G^{a,b}$ has an orientation $D$ with $W(D) \geq M(n)$. 
\end{la}

{\bf Proof:} Denote the sets $\{a_0,a_1,\ldots,a_{n^3}\}$ by $A$ and  
$\{b_0,_1,\ldots,b_{n^3}\}$ by $B$, and the vertex set of $G$ by $V$. \\
To prove the forward implication assume that $G$ has a Hamiltonian path 
$P:v_1, v_2,\ldots,v_n$, where $v_1=a$ and $v_n=b$. 
Orient each edge $v_iv_{i+1}$ of $P$ forward, i.e., as $\overrightarrow{v_iv_{i+1}}$, and orient all 
edges of the form $v_iv_j$ with $i<j-1$ backward, i.e., as $\overrightarrow{v_jv_i}$. In the
resulting orientation we have $d(v_1,v_n)=n-1$. Now orient $a_0a$ and $bb_0$ as $\overrightarrow{a_0a}$ 
and $\overrightarrow{bb_0}$, and orient the edges $a_ia_0$ towards $a_0$ and the edges
$b_ib_0$ towards $b_i$ for $i=1,2,\ldots,n^3$. Denote the resulting digraph by $D$. 
Then $d_D(a_i.b_j) = n+3$,  $d_D(a_0,b_j) = d_D(a_i,b_0)=n+2$ and $d_D(a_0.b_0)=n+1$
for all $i,j \in \{1,2,\ldots,n^3\}$. Hence
\[  W(D) \geq W_D(A,B) = n^6(n+3) + 2n^3(n+2) +n+1  = M(n), \]
as desired. 

For the converse assume that $G$ has order $n$, with $n \geq N_0$ (with $N_0$ to be 
determined later) and that $G^{a,b}$ has an orientation $D$ with $W(D) \geq M(n)$.
We may assume that $\overrightarrow{a_0a} \in D$ since otherwise, if 
$\overrightarrow{aa_0} \in E(D)$ we consider the converse of $D$. 
In order to show that $G$ has Hamiltonian path from $a$ to $b$ it suffices to show that 
\[ d_D(a,b) =n-1. \] 
Suppose not. Then $d_D(a,b) \leq n-2$. We obtain a contradiction by showing that 
this implies that $W(D) <M(n)$. 
Since $\overrightarrow{a_0a} \in E(D)$, there is no path from a vertex in $V\cup B$ to a 
vertex in $A$, so $W_D(V,A)=W_D(B,A)=0$. Hence
\begin{equation} \label{eq:W(D)-as-sum}  
W(D) =  W_D(A\cup B) + W_D(V) + W_D(A,V) + W_D(B,V) + W_D(V,B). 
\end{equation}                   
We first bound $W_D(A \cup B)$. Clearly, $ W_D(A \cup B) = W_D(A) +  W_D(B) + W_D(A,B)$. 
Let $x$ be the number of out-neighbours of $a_0$ in $A$, and 
Let $y$ be the number of in-neighbours of $b_0$ in $B$. 
Then $a_0$ has $n^3-x$ in-neighbours in $A$, and $b_0$ has $n^3-x$ out-neighbours in $B$.
Then $A$ contains $x(n^3-x)$ pairs of vertices at distance $2$ and $n^3$ pairs of vertices
at distance $1$, hence $W_D(A)=2x(n^3-x) + n^3$. Similarly we have 
$W_D(B) = 2y(n^3-y) + n^3$. 
Now $d_D(a_0,b_0) \leq n$ and the distance between the in-neighbours of $a_0$ in $A$ 
and the out-neighbours of $b_0$ in $B$ is at most $n+2$. By a straighforward calculation
we thus get
$W_D(A,B) \leq (n^3+1-x)(n^3+1-y)(n+2) - (2n^3 - x - y +2)$. 
In total we thus obtain 
\[ W_D(A\cup B) \leq  2x(n^3-x) + 2y(n^3-y) + (n^3+1-x)(n^3+1-y)(n+2) + x + y -2, \]
It is easy to see that for $0 \leq x,y \leq n^3$ the right hand side of the above inequality 
is maximised  if $x=y=0$. Hence 
\begin{equation} \label{eq:bound--for-W(AcupB)}
W_D(A \cup B) \leq (n^3+1)^2(n+2) - 2 = n^7 + 2n^6 + O(n^5). 
\end{equation} 
Each of the remaining terms on the right hand side of  \eqref{eq:W(D)-as-sum} is $O(n^5)$. 
Indeed, $W_D(V)$ is a sum of $n(n-1)$ terms, each of which is at most $n-1$. Also 
$W_D(A,V)$ is the sum of $(n^3+1)n$ terms, each of which is not more than $n+1$, so 
$W_D(A,V) = O(n^5)$. Similarly $W_D(B,V) = O(n^5)$ and $W_D(V,B)=O(n^5)$. 
Therefore, \eqref{eq:W(D)-as-sum} implies that 
\begin{equation} \label{eq:bound-for-W(D)}
W(D) \leq n^7 + 2n^6 + O(n^5). 
\end{equation}
Since $M(n) = n^7 + 3n^6 + O(n^5)$, there exists $N_0 \in \mathbb{N}$ such that 
for every $n \in \mathbb{N}$ 
with $n \geq N_0$, the right hand side of \eqref{eq:bound-for-W(D)} is less than $M(n)$.
If $n \geq N_0$, we thus have 
\[ W(D) < M(n), \]
contradicting our assumption that $W(D) \geq M(n)$. The lemma follows. 
\hfill $\Box$

\begin{theo} 
{\sc Wiener-Orientation} is NP-complete.
\end{theo}

{\bf Proof:} 
Our proof is by transformation from {\sc Hamiltonian $(a,b)$-path}. 
Let $N_0$ and $M(n)$ be as in Lemma \ref{la:Hamilton-if-large-W}. 

Given a graph $G$ of order $n$ and vertices $a,b$ of $G$, 

If $n< N_0$, then 
we can determine in constant time if $G$ has a Hamiltonian $(a,b)$-path,
for example by considering all sequences of $n$ vertices if they form an $(a,b)$-path.  

If $n \geq N_0$, then consider $G^{a,b}$. 
Clearly, $G^{a,b}$ can be obtained from $G$ in polynomial time. 
Using $G^{a,b}$ and $M(n)$ as an instance for {\sc Wiener Orientation}, we decide
if $G^{a,b}$ has an orientation of Wiener index at least $M(n)$. By 
Lemma \ref{la:Hamilton-if-large-W}, graph $G$ has a Hammiltonian $(a,b)$-path if and 
only if $G^{a,b}$ has an orientation of Wiener index at least $M(n)$. 
\hfill $\Box$

\section{Orientations of minimum Wiener index}

We conclude this paper by briefly discussing the corresponding minimisation problem,
also raised in \cite{KnoSkrTep2016-2}:
Given a graph $G$, find a (not necessarily strong) orientation of $G$ that minimises
the Wiener index. It was shown by Plesn\'\i k \cite{Ple1984} that this problem becomes
NP-complete if we allow only strong orientations.

We do not know if the problem of finding a (not necessarily strong) orientations of minimum
Wiener index of a given graph is NP-hard.
However, it is easy to see that every orientation of a graph $G$ with $m$ edges has 
Wiener index at least $m$, and below we show that it can be decided in 
polynomial time if a given graph with $m$ edges has an orientation of Wiener index $m$   

A digraph $D$ is transitive if it has the property
that whenever there is a path from a vertex $u$ to a vertex $v$ in $D$, then $D$ contains
the edge $\overrightarrow{uv}$. The following observation is straightforward. 

\begin{obs} \label{obs:transitive-orientation}
Let $D$ be a digraph with $m$ edges. Then $W(D) \geq m$, with equality if and only 
if $D$ is transitive.  
\end{obs}

For a given graph $G$ of order $n$ and size $m$ it follows thus that every orientation 
has Wiener index at least $m$, and that there exists an orientation with Wiener 
index $m$ if and only if $G$ has a transitive orientation. 
It can be decided in time $O(n^2)$ (see for example \cite{Spi1983}) if a given graph has a 
transitive orientation. It follows there exists an algorithm of time complexity 
$O(n^2)$ that decides if $G$ has an orientation of Wiener index $m$.

\end{document}